\documentstyle[amscd]{amsart}
\theoremstyle{plain} \numberwithin{equation}{section}
\newtheorem{thm}{Theorem}[section]
\newtheorem{cor}[thm]{Corollary}
\newtheorem{lem}{Lemma}[section]

\begin{document}
\title[Smooth free involution of $H{\Bbb C}P^3$ and Smith conjecture] {\Large Smooth free involution of
$H{\Bbb C}P^3$ and Smith conjecture for imbeddings of $S^3$ in
$S^6$} \footnote[0]{{\bf Keywords and phrases.}  Smith conjecture,
homotopy complex projective space, involution.
 \endgraf
{\bf 2000 Mathematics Subject Classification.} 57Q45, 55M35,
57R25, 57S17, 57R67.
\endgraf
The first named author is partially supported by the 973  Program of China.
The second named author is partially
 supported by
 the scholar fund of the Ministry of Education in China and
  grants from NSFC (No. 10371020 and No. 10671034).
}
\address{
 Academy of Mathematics and System Science, Academia Sinica, Beijing, 100080,
 P.R. China.}
  \email{libh@@iss06.iss.ac.cn}
 \address{
 Institute of Mathematics, School of Mathematical Science, Fudan University, Shanghai, 200433,
P.R. China.
}
\email{zlu@@fudan.edu.cn}
\author[Bang-he Li and Zhi L\"u]{Bang-he Li and Zhi L\"u  }
 \date{}
\maketitle

{\large
\begin{abstract}
This paper establishes an equivalence between existence of free
involutions on $H{\Bbb C}P^3$ and existence of involutions on
$S^6$ with fixed point set an imbedded $S^3$, then a family of
counterexamples of the Smith conjecture for imbeddings of $S^3$ in
$S^6$ are given by known result on $H{\Bbb C}P^3$. In addition,
this paper also shows that every smooth homotopy complex
projective 3-space admits no orientation preserving smooth free
involution, which answers an open problem \cite{pe}. Moreover, the
study of existence problem for smooth orientation preserving
involutions on $H{\Bbb C}P^3$ is completed.
 \end{abstract}

\section{Introduction}
The original Smith conjecture, which states that no periodic
transformation of $S^3$ can have a tame knotted $S^1$ as its fixed
point set, has been solved in the DIFF category, but it is
generally false in the TOP and PL categories except for some
special cases (see \cite{mb}). However, the generalized Smith
conjecture of codimension two
 is not true in any category (see \cite{gi}, \cite{go},
\cite{su}, and \cite{l1}).
 The generalized Smith conjecture of codimension greater than
two is directly associated with the knot theory of the imbedded
$S^m$ in $S^n$ for $n-m>2$. It is well known that any imbedded
$S^m$ in $S^n$ is unknotted in the TOP and PL categories if
$n-m>2$, and in the DIFF category if $2n>3(m+1)$ (see \cite{ha1},
\cite{le}, \cite{st}, \cite{ze}). Haefliger \cite{ha2}, \cite{ha3}
and Levine \cite{le} showed that there exist infinite imbeddings
of $S^m$ into $S^n$ which are knotted in the DIFF category if
$2n\leq 3(m+1)$ and $m+1\equiv 0 \mod 4$, and that there exists a
knotted $S^{4k+1}$ in $S^{6k+3}$ in the DIFF category. Using
Brieskorn manifolds, some explicit counterexamples for the
generalized conjecture of codimension greater than two were given
in \cite{l2} if $2n\leq 3(m+1)$ and $m+1\equiv 0 \mod 4$ with
$n-m$ being even more than two and $n$ being odd.

\vskip .3cm

The motivation of this paper is to consider
 the generalized Smith conjecture in the DIFF
category for the extreme case $2n=3(m+1)$ with $m+1\equiv 0\mod
4$; especially for $n=6, m=3$. Montgomery and Yang \cite{my}
established the one-to-one
 correspondence $\eta$ between $\Pi$ and $C^3_3$,
 where $\Pi$ is the group
of diffeomorphism classes of all homotopy complex projective
3-spaces, denoted by $H{\Bbb C}P^3$ (for the sum operation, see
\cite{my}), and
 $C^3_3$ is the group of isotopy classes of all imbeddings of
$S^3$ into $S^6$. Note that $C^3_3$ is infinite cyclic (see
\cite{ha3}), so is $\Pi$. This provides a way of dealing with the
Smith conjecture for imbeddings of $S^3$ in $S^6$ by associating
to $H{\Bbb C}P^3$. We will show that there exist infinite distinct
imbeddings $i:S^3\rightarrow S^6$ in $C^3_3$ such that each knot
$(S^6, i(S^3))$ admits a smooth involution, i.e.,  there is a
smooth involution on $S^6$ with $i(S^3)$ as its fixed point set.
This implies that the Smith conjecture for imbeddings of $S^3$ in
$S^6$ is false (see Corollary 1.3 below). Notice that for $k>2$,
there is no ${\Bbb Z}_k$-action on $S^6$ with an imbedded $S^3$ as
the set of fixed points. So the only case concerned with Smith
conjecture is the ${\Bbb Z}_2$-actions. \vskip .3cm Our strategy
to study the ${\Bbb Z}_2$-actions on $S^6$ with an imbedded $S^3$
as fixed point set is to establish its relation with free
involutions on $H{\Bbb C}P^3$, which turns out to be equivalent as
precisely stated by the following theorem.
\begin{thm} \label{t1}
Let $[M]$ be an element in $\Pi$. Then $M$ admits a smooth orientation
reversing free involution
if and only if for any $i:S^3\rightarrow S^6$ in $\eta[M]$, there is a
smooth involution on $S^6$ with $i(S^3)$ as the set of fixed points.
\end{thm}
With respect to $H{\Bbb C}P^3$, Petrie \cite{pe}, Dovermann,
Masuda and Schultz \cite{dms} have already the following
\begin{thm} [\cite{pe}, \cite{dms}] \label{t2}
There are infinitely many homotopy complex projective 3-spaces which
admit a smooth
 orientation reversing free involution.
\end{thm}
Notice that although Petrie's original assertion that every
$H{\Bbb C}P^3$ has an orientation reversing free involution is not
really proved as pointed out in [DMS, p.4], but his proof still
yields Theorem~\ref{t2}.

\vskip .3cm As a consequence of Theorems~\ref{t1} and~\ref{t2}, we
have
\begin{cor}
There exist infinitely many knotted imbeddings of $S^3$ into $S^6$ which offer
counterexamples for the Smith conjecture.
\end{cor}
In addition, we are also concerned with an open problem. In
\cite{pe}, Petrie said that the question of existence of an
orientation preserving
 (free)
 involution on every $H{\Bbb C}P^3$ is still open. The following answers this
question negatively.
\begin{thm} \label{t4}
On every smooth $H{\Bbb C}P^3$, there is no smooth orientation
preserving free involution.
\end{thm}
   Theorem~\ref{t4} means that if an $H{\Bbb C}P^3$ admits a smooth free
   involution, then the involution must be orientation reversing.
   \vskip .3cm
   {\bf Remark.} In his paper \cite{ma}, Masuda studied smooth (nonfree but
   orientation preserving) involutions on $H{\Bbb C}P^3$, and
   proved using Montgomery and Yang  correspondence that every $H{\Bbb C}P^3$
   admits a smooth involution with two copies of ${\Bbb Z}_2$-cohomology
   ${\Bbb C}P^1$ as fixed point set. In contrast to this, he also proved that
   only the standard ${\Bbb C}P^3$ admits a smooth involution with
   ${\Bbb Z}_2$-cohomology ${\Bbb C}P^2$ and ${\Bbb C}P^0$ as fixed point set.
   Thus, Theorem~\ref{t4} with Masuda's results together completes the study of
   existence problem for
    smooth orientation preserving
   involutions on $H{\Bbb C}P^3$.
   \vskip .3cm
The paper is organized as follows. In Section 2, we first give a
proof of Theorem~\ref{t4}, and then prove a basic lemma concerning
 the 2-dimensional homology of the
orbit space of a smooth free involution on an $H{\Bbb C}P^3$. In
Section 3, we show that the surgery processes for Montgomery and
Yang correspondence can still be carried on for ${\Bbb
Z}_2$-actions. This establishes Theorem~\ref{t1}.

\vskip .3cm

{\bf Acknowlegement.} The authors would like to express their
gratitude to M. Masuda for his suggestion of studying the Smith
conjecture by using Montgomery and Yang correspondence, and for
his comments on this paper. The first named author thanks the
Institute of Mathematics of Fudan University for their invitation,
generosity and hospitality, this work is done during his visit to
the institute. The authors also would like to express their
gratitude to the referee, who did an extremely careful reading and
detected some flaws in the original version. The many suggestions
and comments made by him or her considerably improve the
presentation of this paper.

\section{Smooth free involutions of $H{\Bbb C}P^3$}

First, we prove Theorem~\ref{t4} and then give a lemma which is
fundamental for further results.

\vskip .3cm {\em Proof of Theorem~\ref{t4}.} Let $M=H{\Bbb C}P^n$.
Then there is a class $x\in H^2(M;{\Bbb Z})$ such that $H^*(M;{\Bbb
Z})={\Bbb Z}[x]/(x^{n+1})$ and $x^n\in H^{2n}(M;{\Bbb Z})$ is dual
to the fundamental homology class of $M$ (\cite{sp}, Theorem 5, p.
265). If $M$ admits a free involution $\tau$, then $\tau^*(x)=\pm
x$, but $\tau^*(x)=x$ implies that $\tau_*$ is the identity on
$H_*(M;{\Bbb Z})$ and so $\tau$ has a fixed point by the Lefschetz
theorem (\cite{sp}, Theorem 7, p. 195) which is impossible since
$\tau$ is assumed to be free. So $\tau^*(x)=-x$ and $\tau_*$ sends a
generator of $H_2(M;{\Bbb Z})$ into its negative and since
$\pi_2(M)=H_2(M;{\Bbb Z})$, the action of $\pi_1(N)$ on $\pi_2(N)$
is nontrivial (\cite{sp}, Corollary 7, p. 383), where $N$ is the
orbit space of the action $\tau$ on $M$. Since
$\tau^*(x^n)=(-1)^nx^n$, this means that $\tau$ preserves the
orientation if and only if $n$ is even. Thus, when $n=3$, there
exists no orientation preserving free involution on $H{\Bbb C}P^3$.
The proof is completed. \hfill $\Box$

\vskip .3cm

The proof of Theorem~\ref{t4} also gives the following result in
the general case.

\begin{cor}
Suppose that $M$ is an $H{\Bbb C}P^n$ admitting a free involution
$\tau$. Then $\tau$ preserves the orientation if and only if $n$
is even.
\end{cor}

\begin{lem} \label{l1}
Let  $M$ be an $H{\Bbb C}P^3$
 with a smooth free involution $\tau$, and
$p:M\rightarrow N=M/\tau$ be the orbit space projection. Then
$H_2(N;{\Bbb Z})=0.$
\end{lem}

\begin{pf} Let $B$ be a closed M\"obius band with $S^1$ as center, and
$f:S^1\rightarrow N$ be a smooth imbedding representing the
generator of $\pi_1(N)\cong {\Bbb Z}_2$. It is easy to see that
$f$ may extend to a map from $B$ to $N$, still denoted by $f$.
Since $f:\partial B\rightarrow N$ represents twice of the
generator, it is homotopic to zero. Let ${\Bbb R}P^2=B\cup_\lambda
D$, where $D$ is a 2-disk, and $\lambda$ is a diffeomorphism of
$\partial D$ and $\partial B$. Then $f$ extends to a map from
${\Bbb R}P^2$ to $N$, also denoted by
 $f$.
 Change $f$ on the interior
of $D$ by a map $\alpha: S^2\rightarrow N$, we get a map
$f_\alpha:{\Bbb R}P^2\rightarrow N$.

\vskip .3cm

Let ${\Bbb Z}_\xi$ be the local integer coefficient on $N$ twisted
by the line bundle $\xi$ determined by $\tau$. As shown in
\cite{ol}, there is a canonical isomorphism between $H_2(N; {\Bbb
Z}_\xi)$ and $\bar{\Omega}_2(N; \xi)$, so $f$ and $f_\alpha$ may
represent two elements in $H_2(N;{\Bbb Z}_\xi)$, denoted by $[f]$
and $[f_\alpha]$ respectively. Notice that
$\bar{\Omega}_n(X;\phi)$ was defined in \cite{ko} where $\phi$ is
a stable bundle over a space $X$, and Olk in his dissertation
\cite{ol} proved that
$$\bar{\Omega}_n(X;\phi)\cong H_n(X;{\Bbb Z}_\phi) \ \text{for}\
n=0,1,2,3$$ where ${\Bbb Z}_\phi$ is the local integer coefficient
associated to $\phi$ (see also \cite{li}). To see the picture more
clearly, we assume $f:{\Bbb R}P^2\longrightarrow N$ is an
imbedding (this is guaranteed by reason of dimension), and
$f_\alpha=f\sharp \alpha$ is the connected sum of $f$ and an
imbedding $\alpha:S^2\longrightarrow N$ with $\alpha(S^2)\cap
f({\Bbb R}P^2)=\emptyset$. Then $p^{-1}f({\Bbb R}P^2)$ and
$p^{-1}f_\alpha({\Bbb R}P^2)$ are embedded 2-dimensional spheres
in $M$. Let $\alpha':S^2\longrightarrow M$ be such that
$$p\circ\alpha'=\alpha.$$ From the proof of Theorem~\ref{t4}, we know that
 the action of $\pi_1(N)$ on $\pi_2(N)$ is nontrivial, i.e., the action of the generator
of $\pi_1(N)$ sends $x\in{\Bbb Z}$ to $-x\in{\Bbb Z}$. Then the
nontriviality of the action of $\pi_1(N)$ on $\pi_2(N)$ means that
$\alpha'$ and $\tau\circ\alpha'$ represent the elements in
$H_2(M;{\Bbb Z})$ (via Hurewicz homomorphism) with opposite signs.
Geometrically, $p^{-1} f_\alpha({\Bbb R}P^2)$ is the connected sum
of $p^{-1}f({\Bbb R}P^2)$ with $\alpha'$ and $\tau\circ\alpha'$.
When making connected sum with $\alpha'$, the 2-disk $D$ in
$p^{-1}f({\Bbb R}P^2)$ is removed, and when making connceted sum
with $\tau\circ\alpha'$, the 2-disk $\tau(D)$ is removed. Since
$\tau:p^{-1}f({\Bbb R}P^2)\longrightarrow p^{-1}f({\Bbb R}P^2)$ is
orientation reversing, we have
$$[p^{-1}f_\alpha({\Bbb R}P^2)]=[p^{-1}f({\Bbb R}P^2)]+2[\alpha']$$
in $H_2(M;{\Bbb Z})$. By the Gysin homology sequence
\begin{eqnarray} \label{e1}
 0=H_3(M;{\Bbb Z})\buildrel p_*\over\longrightarrow H_3(N;{\Bbb Z})
\longrightarrow H_2(N;{\Bbb Z}_\xi) \buildrel
t\over\longrightarrow
\end{eqnarray}
$$ H_2(M;{\Bbb Z})
\buildrel p_*\over\longrightarrow H_2(N;{\Bbb Z}) \rightarrow
H_1(N;{\Bbb Z}_\xi)
$$where $t$ is the transfer
and the fact that $H_1(N;{\Bbb Z}_\xi)=0$ (which can be seen by
$\bar{\Omega}_1(N;\xi)\cong H_1(N;{\Bbb Z}_\xi)$, and that a map
$g:S^1\longrightarrow N$ with $g^*\xi\cong T(S^1)$ must be
null-homotopic), we see that $H_2(N;{\Bbb Z}_\xi)$ contains at
least one factor of ${\Bbb Z}$, and $H_2(N;{\Bbb Z})=0$ or ${\Bbb
Z}_2$. By another  Gysin homology sequence (\cite{sp}, Problem J,
pp. 282-283)
\begin{eqnarray} \label{e2}
H_2(N;{\Bbb Z})\buildrel t\over\longrightarrow H_2(M;{\Bbb Z})
\buildrel p_*\over\longrightarrow H_2(N;{\Bbb Z}_\xi) \buildrel
\beta \over\longrightarrow
\end{eqnarray}
$$
H_1(N;{\Bbb Z})\cong {\Bbb Z}_2 \longrightarrow H_1(M;{\Bbb Z})=0
$$
we have that $H_2(N;{\Bbb Z}_\xi)$ contains only one factor of
${\Bbb Z}$ and that
$$H_2(N;{\Bbb Z}_\xi)\cong {\Bbb Z}\ \text{or}\ {\Bbb Z}\oplus{\Bbb Z}_2.$$

Next we shall prove that
$$H_2(N;{\Bbb Z}_\xi)\cong{\Bbb Z}\oplus{\Bbb Z}_2$$
is impossible.

\vskip .3cm Suppose that $H_2(N;{\Bbb Z}_\xi)\cong{\Bbb
Z}\oplus{\Bbb Z}_2$. Let $l$ be a generator of $H_2(M;{\Bbb
Z})\cong{\Bbb Z}$. Then there is a pair $(m_0,n_0)$ for
$m_0\in{\Bbb Z}, n_0\in{\Bbb Z}_2$ such that $p_*(l)=(m_0,n_0)$.
If $m_0=0$, then all $(m,0)\in{\Bbb Z}\oplus{\Bbb Z}_2$ with
$m\not=0$ are not in the image of $p_*$, which is impossible by
(\ref{e2}). If $m_0\not=0$, then $(0,1)$ is not in the image of
$p_*$. We claim that \vskip .3cm {\em there exists $(m_1,n_1)$
with $m_1\not=0$ in ${\Bbb Z}\oplus{\Bbb Z}_2$ such that
$(m_1,n_1)$ does not belong to the image of $p_*$ in
$(\ref{e2})$.} \vskip .3cm \noindent Choose an
 $f_\alpha$ defined above having the property
$[p^{-1}f_\alpha({\Bbb R}P^2)]\not=0$, then $[f_\alpha]=(m_1,n_1)$
with $m_1\not =0$. For the homomorphism $\beta$, it can be seen by
the relation of normal bordism groups and homology groups, or by
the Thom isomorphism with local coefficients (\cite{sp}, Problem
J, pp. 282-283) that $\beta(x)=w_1(\xi)\cap x$ for $x\in
H_2(N;{\Bbb Z}_\xi)$ and $w_1(\xi)\in H^1(N;{\Bbb Z}_\xi)$ the
unreduced first Stiefel-Whitney class of $\xi$ as explained in
\cite{ste}. By the first point of view, $\beta[f_\alpha]$ can be
taken as follows. \vskip .2cm For any section $s$ of
$f^*_\alpha\xi$ over ${\Bbb R}P^2$ transversal to the zero
section, let
$$S=\{x\in{\Bbb R}P^2\vert s(x)=0\},$$
then $f_\alpha(S)$ represents $\beta[f_\alpha]$. Obviously we may
take $S$ as being the circle in ${\Bbb R}P^2$ representing the
generator of $\pi_1({\Bbb R}P^2)\cong{\Bbb Z}_2$, thus
$f_\alpha(S)$ represents the generator of $H_1(N;{\Bbb
Z})\cong{\Bbb Z}_2$, and the claim holds. Then by (\ref{e2}) the
image of $\beta$ will have ${\Bbb Z}_2$ as a proper subgroup. This
leads to a contradiction. \vskip .2cm Thus $H_2(N;{\Bbb
Z}_\xi)\cong{\Bbb Z}$, and all $[f_\alpha]$ are odd elements.
\vskip .2cm Now by (\ref{e1}) we have that
\begin{eqnarray} \label{e3}
H_3(N;{\Bbb Z})=0.
\end{eqnarray}
\vskip .2cm If $H_2(N;{\Bbb Z})\cong{\Bbb Z}_2$, then by
(\ref{e3}) and $H_1(N;{\Bbb Z})\cong{\Bbb Z}_2$ and the universal
coefficient theorem, we have that
$$H_0(N;{\Bbb Z}_2)=H_1(N;{\Bbb Z}_2)=H_3(N;{\Bbb Z}_2)={\Bbb Z}_2,\ \ \
H_2(N;{\Bbb Z}_2)={\Bbb Z}_2\oplus{\Bbb Z}_2.
$$
By Poincar\'e duality,
$$H_4(N;{\Bbb Z}_2)={\Bbb Z}_2\oplus{\Bbb Z}_2, \ \ H_5(N;{\Bbb Z}_2)=H_6(N;
{\Bbb Z}_2)={\Bbb Z}_2.$$ Thus the Euler characteristic number of
$N$ would be 3, but it should be half of that of $M$ which is 2.
This contradiction shows that $H_2(N;{\Bbb Z})\cong {\Bbb Z}_2$ is
impossible. \vskip .2cm Thus we have proved that $H_2(N;{\Bbb
Z})=0$.
 This completes the proof.
\end{pf}

\section{Montgomery and Yang correspondence in ${\Bbb Z}_2$-actions}

The proof of Theorem~\ref{t1} is equivalent to establishing the
surgery processes for Montgomery and Yang correspondence under
${\Bbb Z}_2$-actions. For this, we first prove some lemmas.
\begin{lem} \label{l2}
Let $M$ be a smooth $H{\Bbb C}P^3$ with a smooth free involution
$\tau$, and $p:M\rightarrow N=M/\tau$ be the orbit space
projection. Then there is an imbedding $j:{\Bbb R}P^2\rightarrow
N$ such that $p^{-1}j({\Bbb R}P^2)$ is represented by an imbedded
sphere which represents a generator of $\pi_2(M)$.
\end{lem}
\begin{pf} We look at all  $p^{-1}f_\alpha({\Bbb R}P^2)$ stated in the
proof of Lemma~\ref{l1}.  Notice that we have seen in the proof of
Lemma~\ref{l1} that $[f_\alpha({\Bbb R}P^2)]$ are odd elements.
Since $H_2(N;{\Bbb Z})=0$ and $H_3(N;{\Bbb Z})=0$ by
Lemma~\ref{l1} and (\ref{e3}),
 $t$ is an isomorphism in (\ref{e2}),  so
all $[p^{-1}f_\alpha({\Bbb R}P^2)]$ are odd elements in
$H_2(M;{\Bbb Z})=\pi_2(M)={\Bbb Z}$. Thus
 there is
some $\alpha_1$ such that $p^{-1}f_{\alpha_1}({\Bbb R}P^2)$
represents a generator of $\pi_2(M)$. This completes the proof.
\end{pf}

 Recall from \cite{my} that the standard ${\Bbb C}P^3$ can be obtained
by gluing two $S^2\times D^4$ on their boundaries by a map
\begin{eqnarray} \label{e4}
f:S^2\times S^3\rightarrow S^2\times S^3
\end{eqnarray}
where $f(Gu,v)=(Guv,v^{-1})$, $u\longmapsto Gu$ is the Hopf map
$S^3\rightarrow S^2$, and $S^3$ is regarded as the space of unit
quaternions, and $G\subset S^3$ consists of unit complex numbers.
\vskip .3cm
 Regard $S^2$ as the unit sphere in ${\Bbb R}^3={\Bbb R}\times
{\Bbb C}$, and ${\Bbb R}^4$ as the quaternion field which is
identified with ${\Bbb C}^2$ by
$$
x_0+x_1 \text{\bf i}+x_2\text{\bf j}+x_3\text{\bf k}=
x_0+x_1\text{\bf i}+(x_2+x_3\text{\bf i}) \text{\bf
j}\longleftrightarrow (x_0+x_1\text{\bf i}, x_2+x_3\text{\bf
i}).$$ Then the Hopf map is given by ${\Bbb C}^2\rightarrow {\Bbb
R}\times {\Bbb C}$ sending
$$(\psi_1,\psi_2)\longmapsto (\vert\psi_1\vert^2-\vert\psi_2\vert^2,
2\psi_1\bar{\psi}_2).$$ It is easy to check that multiplying
$\text{\bf j}$ on the left side  of the above mapping induces the
antipodal map on $S^2$, i.e., there is a commutative diagram
$$
\begin{CD}
S^3 @>{u\longmapsto \text{\bf j}u}>> S^3 \\
@V{   }VV
@VV{   }V  \\
S^2 @>{Gu\longmapsto -Gu}>> S^2
\end{CD}
$$
Therefore
$$f(G\text{\bf j}u, v)=(G\text{\bf j}uv,v^{-1})=(-Guv,v^{-1})=f(-Gu,v)$$
on $S^2\times D^4$ defines a smooth involution by mapping $(Gu,*)$
to $(-Gu,*)$. This involution commutes with $f$, so there is a
smooth free involution $\tau_0$ on ${\Bbb C}P^3$ such that ${\Bbb
C}P^3/\tau_0$ is glued by two copies of ${\Bbb R}P^2\times D^4$
along their boundaries. Thus, ${\Bbb C}P^3/\tau_0$ is actually a
${\Bbb R}P^2$-bundle over $S^4$. This leads to the following

\begin{lem} \label{l3}
Let $M$ be an $H{\Bbb C}P^3$ with a smooth free involution $\tau$
and $N=M/\tau$. Then  $N_0={\Bbb C}P^3/\tau_0$ and $N$ are
homotopy equivalent.
\end{lem}

\begin{pf}
First we claim that every ${\Bbb R}P^2$-bundle over $S^4$ has a
CW-complex structure such that it contains a cell in each
dimension $i=0, 1, 2, 4, 5, 6$. Every ${\Bbb R}P^2$-bundle over
$S^4$ is the union of two copies of ${\Bbb R}P^2\times D^4$ by
gluing boundaries of two ${\Bbb R}P^2\times D^4$'s. It is
well-known that ${\Bbb R}P^2$ has a CW-decomposition such that it
contains a cell in each dimension $i=0,1,2$, and it is the union
of those three cells. Thus, one of two copies of ${\Bbb
R}P^2\times D^4$ offers one cell in each dimension $i=4,5, 6$. On
the other hand, another of two copies of ${\Bbb R}P^2\times D^4$
has a natural deformation to ${\Bbb R}P^2$, so the $i$-cell for
$i=0,1,2$ is given by this ${\Bbb R}P^2$, and the $i$-cell for
$i=4,5,6$ from the first copy extends to an $i$-cell by the
deformation. Then we obtain the required CW-decomposition.

\vskip .3cm

As a special ${\Bbb R}P^2$-bundle over $S^4$, $N_0$ has such one
CW-decomposition as above. Then we see that the inclusion of
${\Bbb R}P^2$ in $N_0$ induces isomorphisms of homotopy groups up
to dimension 2.  Now, for any $H{\Bbb C}P^3$ space $M$ with free
involution $\tau$, by Lemma~\ref{l2} there is  an imbedding
$j:{\Bbb R}P^2\rightarrow N$ such that the inclusion $j({\Bbb
R}P^2)\hookrightarrow N$ also induces isomorphisms of homotopy
groups up to dimension 2. Let $f$ be a homeomorphism from the
${\Bbb R}P^2$ in $N_0$ to $j({\Bbb R}P^2)$ in $N$. Since the
CW-decomposition of $N_0$ contains no 3-dimensional cells, the
2-skeleton and 3-skeleton of the CW-decomposition of $N_0$ are
equal, and both are just ${\Bbb R}P^2$. Thus, $f$ has been defined
on the boundary of the 4-cell in the CW-decomposition of $N_0$.
Since $\pi_3(N)=0$, the obstruction theory tells us that $f$ can
extend to the 4-cell, denoted still by $f$. Now $f$ is defined on
the 4-skeleton of $N_0$, i.e., on the boundary of the 5-cell of
$N_0$. Since $\pi_4(N)=\pi_5(N)=0$, the same argument as above
shows that finally $f$ can extend to the 6-skeleton of $N_0$,
i.e., $f$ is exactly defined on $N_0$. By the definition of $f$,
$f$ induces isomorphisms of homotopy groups up to dimension 2.
Since $\pi_i(N_0)=\pi_i(N)=0$ for $i=3,4,5, 6$, by  Theorem 3.1 in
page 107 of \cite{hi}, we conclude that $f$ is a homotopy
equivalence.
\end{pf}

\begin{lem} \label{l4}
Let $M$ be an $H{\Bbb C}P^3$ with a smooth free involution $\tau$,
and let $j:{\Bbb R}P^2\rightarrow N=M/\tau$ be the smooth
imbedding described in Lemma~\ref{l2}.  Then the normal bundle of
$j({\Bbb R}P^2)$ is trivial.
\end{lem}
\begin{pf} Since
$${\Bbb C}P^3/\tau_0=({\Bbb R}P^2\times D^4)\cup_\lambda({\Bbb R}P^2
\times D^4)$$ where $\lambda:{\Bbb R}P^2\times S^3\rightarrow
{\Bbb R}P^2\times S^3$ is such that $\lambda:{\Bbb R}P^2\times
v\rightarrow {\Bbb R}P^2\times v^{-1}$ is a homeomorphism, an easy
argument by using Mayer-Vietoris sequence shows that $H_2({\Bbb
C}P^3/\tau_0;{\Bbb Z}_2)={\Bbb Z}_2$ and its generator is
represented by ${\Bbb R}P^2\times *$. The proof of Lemma~\ref{l2}
contains the fact that a generic inclusion $j:j({\Bbb R}P^2)
\longrightarrow N$ induces an isomorphism of homotopy groups at
levels 1 and 2 and so it induces an integral homology isomorphism
at these levels by the Whitehead theorem (\cite{sp}, Theorem 9, p.
399). It follows from the universal coefficient theorems that $j$
induces homology and cohomology isomorphism at levels 1 and 2 with
arbitrary coefficients. In particular, $j({\Bbb R}P^2)$ represents
the generator of $H_2(N;{\Bbb Z}_2)={\Bbb Z}_2$. Since the tangent
bundle $T({\Bbb C}P^3/\tau_0)$ restricted to ${\Bbb R}P^2\times *$
is
$$T({\Bbb C}P^3/\tau_0)\vert_{{\Bbb R}P^2\times *}=
T({\Bbb R}P^2\times *)\oplus\text{trivial bundle}$$ we see that
the Stiefel-Whitney classes $w_1({\Bbb C}P^3/\tau_0)=w_1({\Bbb
R}P^2\times *)\not=0$ and $w_2({\Bbb C}P^3/\tau_0)=w_2({\Bbb
R}P^2\times *)\not=0$. By the proof of Lemma~\ref{l3},  a homotopy
equivalence $f: {\Bbb C}P^3/\tau_0 \longrightarrow M/\tau=N$ can
be chosen so that $f:{\Bbb R}P^2\times *\longrightarrow j({\Bbb
R}P^2)$ is a homeomorphism.

\vskip .3cm

It is well known that Stiefel-Whitney classes are homotopy
invariant, so $w_1(N)\not=0$ in $H^1(N;{\Bbb Z}_2)$ and
$w_2(N)\not=0$ in $H^2(N;{\Bbb Z}_2)={\Bbb Z}_2$. Thus
$w(T(N)\vert_{j({\Bbb R}P^2)})=w(j({\Bbb R}P^2))$, and the
Stiefel-Whitney classes of the normal bundle are trivial. This
means that the normal bundle of $j({\Bbb R}P^2)$ is stably trivial
since its Stiefel-Whitney classes vanish at levels 1 and 2 and
$j({\Bbb R}P^2)$ is a 2-complex. Thus, the normal bundle of
$j({\Bbb R}P^2)$ is trivial since its fiber has dimension 4 (which
is strictly larger than 2).  This completes the proof.
\end{pf}

With the above understood, we are going to complete the proof of
Theorem~\ref{t1}. \vskip .3cm

{\em Proof of Theorem~\ref{t1}.} Suppose that $M$ admits  a smooth
free involution $\tau$. Then by Lemma~\ref{l4} there is an
imbedding $k:{\Bbb R}P^2\times D^4\rightarrow M/\tau$. Hence there
is an imbedding $h:S^2\times D^4\rightarrow M$ such that
$$h(S^2\times D^4)=p^{-1}k({\Bbb R}P^2\times D^4)$$
where $p:M\rightarrow M/\tau$ is the projection.  By
Lemma~\ref{l2}, $h:S^2\times 0\rightarrow M$ is a primary
imbedding defined in \cite{my}. Now, by Lemma 5 in \cite{my},
there is another primary imbedding $h':S^2\times D^4\rightarrow M$
such that
$$h(S^2\times D^4)\cup h'(S^2\times D^4)=M$$
and
$$h(S^2\times D^4)\cap h'(S^2\times D^4)=h(S^2\times S^3)=
h'(S^2\times S^3).$$
Clearly, on $h(S^2\times D^4)$,
$$\tau h(x,y)=h(-x,y).$$
 According to [MY, Appendix and Lemma 11], an imbedding
$i:S^3\rightarrow S^6$ is given by
$$S^6=(M-\text{int}h(S^2\times D^4))\cup_\alpha(D^3\times S^3)$$
$$i:S^3=\{0\}\times S^3\subset D^3\times S^3\subset S^6$$
where $\alpha: h(S^2\times S^3)\rightarrow \partial(D^3\times S^3)=S^2\times
 S^3$ is defined by $\alpha=fh^{-1}$ and $f(Gu,v)=(Guv,v^{-1})$.
 \vskip .2cm
 Now on $(M-\text{int}h(S^2\times D^4))=h'(S^2\times D^4)\subset S^6$,
 there is
 the involution $\tau_1(=\tau)$.
We define an involution $\tau_2$ on $D^3\times S^3$ by $\tau_2(x,y)=(-x,y)$.
Since on $h(S^2\times S^3)$, we have
$$\alpha\tau_1 h(x,y)=fh^{-1}\tau h(x,y)=fh^{-1}h(-x,y)=f(-x,y)$$
and
$$\tau_2\alpha h(x,y)=\tau_2 fh^{-1}h(x,y)=\tau_2f(x,y).$$
Represent $(x,y)\in S^2\times S^3$ by $(Gu,v)$, then
$$f(Gu,v)=(Guv,v^{-1})$$
and
$$\tau_2 f(Gu,v)=(-Guv,v^{-1})=f(-Gu,v).$$
Hence
$$\alpha\tau_1=\tau_2\alpha$$
 on $h(S^2\times S^3)$. Combining
$\tau_1$ with $\tau_2$ together, we obtain then an involution on $S^6$ with
$i(S^3)$ as the set of fixed points.
\vskip .3cm
Conversely, let $i:S^3\rightarrow S^6$ be an imbedding such that there is
an involution $\tau$ on $S^6$ with $i(S^3)$ as the set of fixed points.
Take a $\tau$-invariant Riemannian metric on $S^6$,
then $\tau$ induces a
bundle map of the normal bundle of $i(S^3)$,
i.e., the orthogonal
bundle of $T(i(S^3))$ in $T(S^6)$ according to the Riemannian metric,
 covering the identity map of
$i(S^3)$, and on every fiber of the normal bundle of $i(S^3)$, the
bundle map must set $v$ to $-v$. Then by using the exponential
map, it follows that there is an ${\Bbb Z}_2$-equivariant
imbedding $h:D^3\times S^3\rightarrow S^6$ such that
$$h(0,x)=i(x)\ \text{ and }\ \tau h(v,x)=h(-v,x).$$
By \cite{my}, there is an imbedding $k:S^2\times D^4\rightarrow
S^6$ such that
$$S^6=k(S^2\times D^4)\cup h(D^3\times S^3)$$
and
$$k(S^2\times D^4)\cap h(D^3\times S^3)=k(S^2\times S^3)=h(S^2\times S^3).$$
By the proof of Lemma 12 in \cite{my}, we may change $h$ by a map
$\mu:S^3\rightarrow SO(3)$ to get an imbedding $k':D^3\times
S^3\rightarrow S^6$, so that $k$ and $k'$ satisfy the condition of
[MY, Lemma 12], and
$$\tau k'(v,x)=k'(-v,x).$$
 Now let $\lambda: k(S^2\times S^3)\rightarrow S^2\times S^3$ be the map
 defined by $\lambda=f{k'}^{-1}$, where $f$ is the map stated in (\ref{e4}), then
 $$M=k(S^2\times D^4)\cup_\lambda(S^2\times D^4)$$
 is an $H{\Bbb C}P^3$ by [MY, Lemma 12], and $\eta[M]=[i]$.
 \vskip .2cm
 To see that $M$ has a free involution, let $\tau_1=\tau$ on $k(S^2\times D^4)$, and $\tau_2$ on $S^2\times D^4$ is given by
 $$\tau_2(v,x)=(-v,x).$$
 Then on $k(S^2\times S^3)=k'(S^2\times S^3)$,
$\lambda \tau_1 k'(v,x)=\lambda k'(-v,x)=f(-v,x)$
 and
 $$\tau_2\lambda k'(v,x)=\tau_2 f(v,x).$$
 As we have seen before,
 $$\tau_2 f(v,x)=f(-v,x)$$
  so
  $$\lambda\tau_1
 =\tau_2\lambda$$
  on $k(S^2\times S^3)$. Therefore,
 combining $\tau_1$ and $\tau_2$ gives a free involution on $M$.
 This completes the proof.
 \hfill $\Box$

\end{document}